\documentclass[12pt]{amsart}

\usepackage{graphicx}
\usepackage{amsmath}
\usepackage{amssymb}
\usepackage{amsfonts}
\usepackage{verbatim}
\usepackage{scalefnt}
\usepackage{multirow}
\usepackage{enumerate}
\usepackage{mathtools}
\usepackage{float}
\usepackage{enumitem}
\restylefloat{figure}
\usepackage[margin=3cm]{geometry}

\usepackage{fancyhdr}
\usepackage{hyperref}

\usepackage{bbm}

\newtheorem{problem}{Problem}
\newtheorem{theorem}{Theorem}
\newtheorem{lemma}[theorem]{Lemma}

\newtheorem{corollary}[theorem]{Corollary}
\newtheorem{conjecture}[theorem]{Conjecture}

\newtheorem*{definition*}{Definition}

\newtheorem{remark}[theorem]{Remark}

\numberwithin{equation}{section}
\numberwithin{theorem}{section}

\title[Equidistribution of  subset sums]%
  {Equidistribution of  subset sums}

\author{P\'eter P\'al Pach}
\email{pachpp@renyi.hu}
\address{HUN-REN Alfr\'ed R\'enyi Institute of Mathematics, Re\'altanoda utca 13--15., H-1053 Budapest,  Hungary; MTA--HUN-REN RI Lend\"ulet ``Momentum'' Arithmetic Combinatorics Research Group, Re\'altanoda utca 13--15., H-1053 Budapest,  Hungary; Department of Computer Science and Information Theory, Budapest University of Technology and Economics, M\H{u}egyetem rkp. 3., H-1111 Budapest, Hungary.}

\thanks{}

\begin{document}

\begin{abstract}	
We answer a question of Katona and Makar-Limanov, by showing that in an abelian group of order $2h$ the $h$-element subset sums are asymptotically (as $h\to \infty$) equidistributed. In fact we prove a more general result where the order of the group can be arbitrary, also providing a bound for the ``error term''.
\end{abstract}

\date{\today}
\maketitle

\section{Introduction} 

In this note, we address a question of Katona and Makar-Limanov~\cite{KML08} about the distribution of subset sums in abelian groups. Let $G$ be an abelian group of order $n$.  Let us set $h:=\left[ \frac{n}{2}\right]$ and define the family
$$\mathcal{F}_a=\{ \{x_1,\dots,x_h\}:\ x_1,\dots,x_h\in G\text { are distinct and }x_1+\dots+x_h=a\}.$$
Observe that $\sum \limits_{a\in G} |\mathcal{F}_a|=\binom{n}{h}$. Katona and Makar-Limanov asked how small and large $|\mathcal{F}_a|$ can be. Let $$m(G):=\min\limits_{a\in G} |\mathcal{F}_a|\text{ and }M(G):=\max\limits_{a\in G} |\mathcal{F}_a|,$$
also $M(n):=\max\limits_{|G|=n}M(G)$. Katona and Makar-Limanov showed that 
$$\frac{1}{n}\binom{n}{h}\leq M(G)\leq \frac{2}{n}\binom{n}{h}(1+o(1))$$
and formulated the following problems:
\begin{problem}
    Determine 
    $$\limsup \frac{M(n)}{\frac{1}{n}\binom{n}{h}}.$$
\end{problem}
\begin{problem}
Is
$$\lim\limits_{|G|\to\infty} \frac{m(G)}{M(G)}=1$$
always true?
\end{problem}
\begin{problem}
How small can $m(G)$ be in asymptotical sense?
\end{problem}
They also mentioned that the following conjecture is widely  believed to be true.
\begin{conjecture} (folklore)
$$\lim\limits_{n\to \infty}
\frac{m(\mathbb{Z}_n)}{
M(\mathbb{Z}_n)}
= 1$$
\end{conjecture}
Their motivation was coming from a coding theoretic problem, which asks for the maximum size of a binary code of length $n$, fixed weight $h$ with minimum Hamming-distance 4. Observe that the set system $\mathcal{F}_a$ can be used to construct such a code, since the symmetric difference of any two sets in $\mathcal{F}_a$ is at least 4.  For some applications of this coding problem to combinatorics, see \cite{BK07,GS80}
and the survey \cite{Kat08}.

Katona and Makar-Limanov~\cite{KML08} resolved the special case $G=\mathbb{Z}_2^r$.

We solve the above-mentioned problems for arbitrary abelian groups by showing that $$|\mathcal{F}_a|=(1+o(1))\frac{1}{n}\binom{n}{h}$$
(as $n\to \infty$) for every $a\in G$. 
 
Let us note that the interesting case is when $n=2h$ is even. Assume that $n=2h+1$ is odd. Observe that by shifting all elements in an $h$-tuple $\{x_1,\dots,x_h\}\in \mathcal{F}_a$ by $d$ we get the $h$-tuple $\{x_1+d,\dots,x_h+d\}$ which is contained in $\mathcal{F}_{a+hd}$. This yields that $|\mathcal{F}_a|=|\mathcal{F}_{a+hd}|$ for every $a,d\in G$. In particular, $$|\mathcal{F}_0|=|\mathcal{F}_{hd}|=|\mathcal{F}_{2hd}|=|\mathcal{F}_{-d}|$$
for every $d\in G$, implying that all the sets $\mathcal{F}_a$ have the exact same size $\frac{1}{n}\binom{n}{h}$. 

While the original coding theoretic motivation suggested investigating the cases $n=2h$ and $n=2h+1$, it is natural to look at other instances, as well. We may restrict our attention to the cases $h\leq n/2$, since an equidistribution result for $h=h_0$ implies the analogous result for $h=n-h_0$, too. The case  $h=1$ is trivial, and for $h=2$ the sums $x_1+x_2$ (with distinct $x_1,x_2$) may not be equidistributed asymptotically. For instance, for $G=\mathbb{Z}_2^r$ the element 0 is not represented at all in this form. Our methods provide asymptotic equidistribution results for the whole range $3\leq h\leq n-3$ for arbitrary abelian groups.

It will be more convenient to work with $h$-tuples $(x_1,\dots,x_h)$ instead of subsets $\{x_1,\dots,x_h\}$. Each $h$-element subset $\{a_1,\dots,a_h\}\in \mathcal{F}_{a}$ corresponds to $h!$ solutions of the equation 
$$x_1+\dots+x_h=a$$
consisting of distinct elements: we can take all the $h!$ permutations of $a_1,\dots,a_h$.

We prove the following results.

\begin{theorem}\label{thm-small-h}
Let $h\geq 3$ be a fixed integer. Let $G$ be an abelian group of order $n\geq 2h$. Then the number of solutions to $x_1+\dots+x_h=a$ with distinct $x_1,\dots,x_h\in G$ is $\frac{1}{n}\cdot \frac{n!}{(n-h)!}+O_h(n^{h/2})$. Therefore, the $h$-element subset sums are asymptotically equidistributed in $G$ (as $n\to\infty$). 
\end{theorem}

In our second theorem $h$ can also grow.

\begin{theorem}\label{thm-main}
Let $h$ be a positive integer and $G$ be an abelian group of order $n$. If $\min(h,n-h)$ is sufficiently large, then the number of solutions to $x_1+\dots+x_h=a$ with distinct $x_1,\dots,x_h\in G$ is $\frac1n\cdot\frac{n!}{(n-h)!}+\Delta$, where $|\Delta|\leq \left(\frac34 \right)^h \cdot \frac1n\cdot\frac{n!}{(n-h)!} $.

Therefore, the $h$-element subset sums are asymptotically equidistributed in $G$ (as $h,n\to \infty$).
\end{theorem}
This readily implies the following statement resolving Problems 1-3:
\begin{corollary}
  Let $G$ be an abelian group of order $n=2h$. For every $a\in G$ the number of $h$-element sets $\{x_1,\dots,x_h\}\subseteq G$ such that $x_1+\dots+x_h=a$ is $(1+O(0.75^h))\frac{1}{n}\binom{n}{h}$.
\end{corollary}

Hence, the $\limsup$ in Problem 1 is equal to 1, the answer to Problem 2 is affirmative, $m(G)\geq (1-o(1))\frac{1}{n}\binom{n}{h}$ and Conjecture 4 holds for arbitrary abelian groups.

The author learned -- through Gyula O.~H.~Karona -- that Jing Wang  in a parallel work~\cite{Wan25} also obtained similar results solving Problems 1-3. As our techniques are different, we decided to keep our results in two separate papers.

In this paper the standard notation  $O$ is applied to positive quantities in the usual way.
That is,  $Y = O(X)$ means that $X \geq cY$,
for some absolute constant $c > 0$, we also use the Landau $o$ notation.

\section{Preliminary lemmas}

In this section we present two lemmas that we will use in our proof.

Note that the number of solutions to 
\begin{equation}\label{eq}
x_1+\dots+x_h=a\ (x_1,\dots,x_h\in G)
\end{equation}
{\it without} the distinctness condition is $n^{h-1}$ for every $a$, as $x_1,\dots,x_{h-1}$ can be chosen arbitrarily, then $x_h=a-(x_1+\dots+x_{h-1})$  is uniquely determined by \eqref{eq}. In order to be able to count solutions where $x_1,\dots,x_h$ are distinct we prove an inclusion-exclusion formula for Kronecker-delta functions. 

Let us introduce some notation. Let $\mathcal P$ be a partition of the indices $1,\dots,h$. If $i$ and $j$ belong to the same equivalence class, we write $i \sim_{\mathcal P} j$. Let $k_1,\dots,k_r$ denote the sizes of the equivalence classes in $\mathcal P$. Define 
$$c(\mathcal P):=(-1)^{\sum\limits_{i=1}^r (k_i-1)}\prod\limits_{i=1}^r (k_i-1)!,$$
in particular $c(\mathcal P)=1$, if $1,\dots,h$ are all in 1-element classes.

Let $\delta(\mathcal P)(x_1,\dots,x_h):=1$, if $x_i=x_j$ for every $i$ and $j$ such that $i\sim_{\mathcal P} j$ and $\delta(\mathcal P):=0$ otherwise. We say that the $h$-tuple $(a_1,\dots,a_h)$ is compatible with $\mathcal P$ if $\delta(\mathcal P)(a_1,\dots,a_h)=1$. For instance, if $\mathcal P$ is the partition $\{1,2\},\{3,4,5\},\{6\}$, then $(a_1,a_2,a_3,a_4,a_5,a_6)$ is compatible with $\mathcal P$, if $a_1=a_2$ and $a_3=a_4=a_5$ (note that the values $a_1,a_3,a_6$ may not be distinct).

\begin{lemma}\label{lem-incl-excl}
The number of solutions to $x_1+\dots+x_h=a$ with distinct $x_1,\dots,x_h\in G$ is
$$\sum\limits_{\mathcal P}\sum\limits_{x_1,\dots,x_h\in G} c(\mathcal{P})\delta(\mathcal{P}) \mathbbm{1}_{x_1+\dots+x_h=a}(x_1,\dots,x_h),$$
where the summation is taken for all partitions of $1,2,\dots,h$ and $\mathbbm{1}_{x_1+\dots+x_h=a}(x_1,\dots,x_h)=1$, if $x_1+\dots+x_h=a$ and it is 0, otherwise.
\end{lemma}

\begin{remark}
 Let us illustrate the statement by the case $h=4$. In this case there are 15 partitions. 
 \begin{itemize}
    \item One partition, where we have four singletons: $\{1\},\{2\},\{3\},\{4\}$.
    \item Six partitions with one pair and two singletons, e.g. let us denote $\{1,2\},\{3\},\{4\}$ by $12$.
    \item Three partitions consisting of two pairs, e.g. let us denote $\{1,2\},\{3,4\}$ by $12,34$.
    \item Four partitions consisting of a triplet and a singleton, e.g. let us denote $\{1,2,3\},\{4\}$ by $123$.
    \item One partition consisting of one set $\{1,2,3,4\}$, denoted by $1234$.
 \end{itemize}
 Then the value of the function 
 $$1-\delta_{12}-\delta_{13}-\delta_{14}-\delta_{23}-\delta_{24}-\delta_{34}+\delta_{12,34}+\delta_{13,24}+\delta_{14,23}+2\delta_{123}+2\delta_{124}+2\delta_{134}+2\delta_{234}-6\delta_{1234}$$
 is 1, if $x_1,x_2,x_3,x_4$ are all distinct and 0 otherwise.

 We did not find this lemma in the literature, so we provide a proof. For applications of Kronecker-delta functions for dealing with distinctness, see also \cite{BMPS,Nas20}.
\end{remark}

\begin{proof} [Proof of Lemma~\ref{lem-incl-excl}]

If $a_1,\dots,a_h$ are all distinct, then $\delta(\mathcal P)=0$ except the case when all the $h$ indices are in 1-element classes. Thus the solutions $(a_1,\dots,a_h)$ consisting of distinct elements are indeed counted exactly once.

Now, we show that otherwise the total count for $(a_1,\dots,a_h)$ is 0. Let us take an $h$-tuple $(a_1,\dots,a_h)$ satisfying $a_1+\dots+a_h=a$ and denote by $\mathcal P$ the partition where $i\sim_{\mathcal P}j$ if and only if $a_i=a_j$. We shall prove that $\sum\limits_{\mathcal P'\leq \mathcal P} c(\mathcal P')=0$, since, if $\mathcal P'$ is not a refinement of $\mathcal P$, then $\delta(\mathcal P')(a_1,\dots,a_h)=0$. Let $k_1,\dots,k_r$ be the sizes of the equivalence classes larger than 1 in $\mathcal P$. We prove the statement by induction on $k_1+\dots+k_r$. 

The base case $k_1+\dots+k_r=2$ is trivial ($r=1,k_1=2$ and the sum is $1-1=0$), so let us prove the induction step.

First, assume that $r>1$. Let $I_1,\dots,I_r$ be the equivalence classes (containing at least two elements) in $\mathcal P$. Each refinement $\mathcal P'\leq \mathcal P$ can be expressed as $\bigcup \mathcal P_i'$, where $\mathcal P_i'$ is a partition of the elements in $I_i$. Let $\mathcal P_i$ denote the partition of $I_i$ where there is only one equivalence class of size $|I_i|$. Observe that $c(\mathcal P')=\prod\limits_{i=1}^r c(\mathcal P_i')$, thus $\sum\limits_{\mathcal P'\leq \mathcal P} c(\mathcal P')=\prod\limits_{i=1}^r \sum\limits_{\mathcal P_i'\leq \mathcal P_i} c(\mathcal P_i')=0$, by the induction hypothesis applied to $\mathcal P_i$ for some $1\leq i\leq r$.

Now, assume that $r=1$, for brevity let $k:=k_1$. The smallest possible value for $k$ is $k=2$, which was our base case, so we may assume that $h=k\geq 3$. Let $I_h$ denote the equivalence class of $h$. For a fixed set $I=I_h$ the remaining indices in $[h]\setminus I=:J$ can be partitioned arbitrarily, and \begin{equation}\label{eq-part}\sum\limits_{\mathcal P:\ I_h=I} c(\mathcal P)=(-1)^{|I|-1}(|I|-1)!\sum\limits_{\mathcal P' \text{ is a partition of $J$}} c(\mathcal P').
\end{equation}
By the induction hypothesis, this sum is 0 unless $J=\emptyset$ or $|J|=1$.

If $J=\emptyset$, then the sum in \eqref{eq-part} is $(-1)^{h-1}(h-1)!$, while for $|J|=1$ each of the $h-1$ summands are equal to $(-1)^{h-2}(h-2)!$, yielding that the total sum is indeed 0. This concludes the proof of the lemma. 

\end{proof}

Another key observation is that in Lemma~\ref{lem-incl-excl} the contribution from {\it most} partitions $\mathcal P$ is independent of the value of $a$.  

\begin{lemma}\label{lem-equi}
Let $k_1,\dots,k_r$ be the sizes of the equivalence classes (including 1-element classes) in the partition $\mathcal P$ of the elements $1,2,\dots,h$. If 
$$\gcd(k_1,\dots,k_r,n)=1,$$
then
$$\sum\limits_{x_1,\dots,x_h\in G} \delta(\mathcal P) \mathbbm{1}_{x_1+\dots+x_h=a}(x_1,\dots,x_h)=n^{r-1}.$$
\end{lemma}

\begin{proof}
Let $S(a)$ denote the number of those solutions to $x_1+\dots+x_h=a$ that are compatible with $\mathcal P$:
$$S(a):=\sum\limits_{x_1,\dots,x_h\in G} \delta(\mathcal P) \mathbbm{1}_{x_1+\dots+x_h=a}(x_1,\dots,x_h).$$
Note that $S(a)$ is the number of $r$-tuples $~{(y_1,\dots,y_r)\in G^r}$ such that  
$$k_1y_1+\dots+k_ry_r=a.$$
Notice that the multiset $\{k_iy: y\in G\}$ contains only elements from the subgroup $k_iG=\{k_ig:g\in G\}$, each of them is contained with the same multiplicity $|G|/|k_iG|$. Therefore, $S(a)=S(a')$ if $a-a'\in k_iG$. Since the condition $\gcd(k_1,\dots,k_r,n)=1$ implies that  $$k_1G+k_2G+\dots+k_rG=G,$$ 
the sums $k_1y_1+\dots+k_ry_r$ are indeed equidistributed, thus $S(a)=n^r/n=n^{r-1}$. 

\end{proof}

\section{Proofs}

Now, we are ready to present the proofs of Theorem~\ref{thm-small-h} and Theorem~\ref{thm-main}.

\begin{proof}[Proof of Theorem~\ref{thm-small-h}]
According to Lemma~\ref{lem-incl-excl} the number of solutions $x_1+\dots+x_h=a$ with distint $x_1,\dots,x_h\in G$ is
$$T(a):=\sum\limits_{\mathcal P}\sum\limits_{x_1,\dots,x_h\in G} c(\mathcal{P})\delta(\mathcal{P}) \mathbbm{1}_{x_1+\dots+x_h=a}(x_1,\dots,x_h).$$
Lemma~\ref{lem-equi} yields that the contribution from those partitions $\mathcal{P}$ where the greatest common divisor of the sizes of the equivalence classes is relatively prime to $n$ is the same for every $a\in G$. 

For the remaining partitions $\mathcal{P}$ there exists some $2\leq d\mid \gcd(h,n)$ such that each equivalence class has size divisible by $d$, so the number of equivalence classes is at most $h/d\leq h/2$. Hence, the contribution of such a partition $\mathcal P$ to $T(a)$ is at most $|c(\mathcal P)|\cdot n^{h/2}$, thus the total contribution of these partitions is at most $c_h\cdot n^{h/2}$ for some constant $c_h$ depending only on $h$. (For instance, the choice $c_h=\sum\limits_{\mathcal P} |c(\mathcal P)|$ is fine.)

For $h>2$ the error term $O(n^{h/2})$ is negligible compared to $\frac{1}{n}\cdot\frac{n!}{h!}$.
\end{proof}

\begin{proof}[Proof of Theorem~\ref{thm-main}]
According to Lemma~\ref{lem-incl-excl} the number of solutions $x_1+\dots+x_h=a$ with distint $x_1,\dots,x_h\in G$ is
$$T(a):=\sum\limits_{\mathcal P}\sum\limits_{x_1,\dots,x_h\in G} c(\mathcal{P})\delta(\mathcal{P}) \mathbbm{1}_{x_1+\dots+x_h=a}(x_1,\dots,x_h).$$
Lemma~\ref{lem-equi} yields that the contribution from those partitions $\mathcal{P}$ where the greatest common divisor of the sizes of the equivalence classes is relatively prime to $n$ is the same for every $a\in G$.  (Note that this observation already provides a complete answer to the previously discussed easy case $n=2h+1$.)

Let us bound the contribution coming from the remaining partitions $\mathcal{P}$. In case of these there exists some $2\leq d\mid \gcd(h,n)$ such that each equivalence class has size divisible by $d$. In such a partition let $\alpha_i$ denote the number of classes of size $i d$. Note that $\sum\limits_{i=1}^h i\alpha_i=h/d$. 

Observe that $|c(\mathcal P)|=\prod\limits_{i=1}^h [(id-1)!]^{\alpha_i}$ and the total number of those $h$-tuples $(x_1,\dots,x_h)\in G$ that are compatible with $\mathcal P$ is $n^{\sum \alpha_i}$. The total number of these partitions is 
$$\frac{h!}{\prod\limits_{i=1}^h \left[(\alpha_i)![(id)!]^{\alpha_i}\right]},$$
thus the total contribution from this type of partitions to $T(a)$ is at most
$$\frac{h!}{\prod\limits_{i=1}^h \left[(\alpha_i)![(id)!]^{\alpha_i}\right]}\cdot \prod\limits_{i=1}^h [(id-1)!]^{\alpha_i}\cdot n^{\sum\alpha_i}=h!\cdot\prod\limits_{i=1}^h\frac{n^{\alpha_i}}{ (\alpha_i)!(id)^{\alpha_i}}
$$
in absolute value. 

Now, we want to bound the sum of the products $\prod\limits_{i=1}^h\frac{n^{\alpha_i}}{ (\alpha_i)!(id)^{\alpha_i}}$ to all nonnegative integers $\alpha_i$ ($i\geq 1$) satisfying $\sum\limits_{i=1}^h i\alpha_i=h/d$. By using the estimate $(t/e)^t\leq t!$ we have 
$$\frac{n^{\alpha_i}}{ (\alpha_i)!(id)^{\alpha_i}}\leq \frac{n^{\alpha_i}}{ (\alpha_i/e)^{\alpha_i}(id)^{\alpha_i}}= \frac{(ne)^{\alpha_i}}{(di\alpha_i)^{\alpha_i}} .$$

\noindent
(For $\alpha_i=0$ we define $0!=0^0=1$, as usual.) 

We will show a bound of the form 
\begin{equation}\label{eq-1.2}
    \prod\limits_{i=1}^h \frac{(\lambda_ine)^{\alpha_i}}{(di\alpha_i)^{\alpha_i}} \leq M,
\end{equation}
if $h$ is sufficiently large with $1<\lambda_i:=e^{i/1000}$ and $M:=\left(\frac34\right)^h\frac{1}{Chn}\binom{n}{h}$, where $C:=\prod\limits_{i=1}^\infty \frac{\lambda_i}{\lambda_i-1}<\infty$.

Note that this is equivalent to
$$\prod\limits_{i=1}^h \frac{(ne)^{\alpha_i}}{(di\alpha_i)^{\alpha_i}} \leq  \left(\prod\limits_{i=1}^h \lambda_i^{-\alpha_i}\right)M.$$
Instead of summing this for each $h$-tuple satisfying $\sum\limits_{i=1}^h i\alpha_i=h/d$, we may take the sum for every $h$-tuple $(\alpha_1,\dots,\alpha_h)$ such that the  $\alpha_i$ are nonnegative integers, getting the bound
\begin{multline}\label{eq-sum}
\sum\limits_{\substack{\alpha_1,\dots,\alpha_h\in \mathbb{Z}_{\geq 0}:\\ \alpha_1+\dots+\alpha_h=h/d}}\prod\limits_{i=1}^h \frac{(ne)^{\alpha_i}}{(di\alpha_i)^{\alpha_i}} \leq\sum\limits_{\alpha_1,\dots,\alpha_h\in \mathbb{Z}_{\geq 0}} 
  \prod\limits_{i=1}^h  \frac{(ne)^{\alpha_i}}{(di\alpha_i)^{\alpha_i}} \leq \sum\limits_{\alpha_1,\dots,\alpha_h\in \mathbb{Z}_{\geq 0}} \left(\prod\limits_{i=1}^h \lambda_i^{-\alpha_i}\right)M\leq  \\ 
 \leq \left(\prod\limits_{i=1}^h \sum\limits_{\alpha_i=0}^{\infty} \lambda_i^{-\alpha_i}\right) M 
  = \left(\prod\limits_{i=1}^h \frac{\lambda_i}{\lambda_i-1} \right)M=CM.
\end{multline}

As $d\mid h$, adding up these for every choice of $d$ we get the upper bound 
\begin{equation}\label{eq-final}
\left|T(a)-\frac1n\cdot\frac{n!}{(n-h)!}\right|\leq h!\cdot hCM=\left(\frac34\right)^h\cdot \frac{1}{n}\cdot\frac{n!}{(n-h)!}
\end{equation}
for every $a\in G$.

Now, let us prove \eqref{eq-1.2}. 
Clearly, we may assume that $d=2$:
$$\prod\limits_{i=1}^h \frac{(\lambda_ine)^{\alpha_i}}{(2i\alpha_i)^{\alpha_i}} \leq M.$$
Setting $\beta_i:=i\alpha_i/h$ we have $\sum \beta_i=1/2$ and want to estimate the product
$$\prod\limits_{i=1}^h \frac{(\lambda_ine)^{\beta_ih/i}}{(2\beta_ih)^{\beta_ih/i}} .$$

Let us estimate the logarithm of this product. (For $\beta_i=0$ the expression $\beta_i \log \beta_i$ should be considered 0 in the following calculation.)

Observe that
\begin{equation}\label{eq-1}
\sum\limits_{i=1}^h \frac{\beta_ih}{i}\log(n/h)\leq \frac{h}{2} \log(n/h),
\end{equation}
since $\sum\limits_{i=1}^h \frac{\beta_i}{i}\leq \sum\limits_{i=1}^h \beta_i=1/2$.

Also,
\begin{equation}\label{eq-2}\sum\limits_{i=1}^h \frac{\beta_i h}{i} (\log\lambda_i+1-\log2)\leq \frac{1-\log2}{ 2}h+
h\sum\limits_{i=1}^h \frac{\beta_i\log \lambda_i}{i}\leq 0.155h,
\end{equation}
as we set $\lambda_i=e^{i/1000}$.

Finally,
\begin{equation}\label{eq-3}
h\sum\limits_{i=1}^h \frac{-\beta_i \log \beta_i}{i}
\end{equation}
is maximized when $\frac{1+\log \beta_i}{i}$ is the same for every $i$, yielding the bound
$$h\sum\limits_{i=1}^h \frac{-\beta_i\log\beta_i}{i}\leq \left(\frac12\log\frac{e+2}{e}+\frac{1}{e}\log \frac{e+2}{2}\right)h\leq 0.592h.$$

Hence, from \eqref{eq-1}, \eqref{eq-2} and \eqref{eq-3} we obtain that
$$\prod\limits_{i=1}^h \frac{(\lambda_ine)^{\beta_ih/i}}{(2\beta_ih)^{\beta_ih/i}} \leq e^{3h/4}(n/h)^{h/2}.$$

Next, we show that 
\begin{equation}\label{eq-4}
e^{3h/4}\left(\frac{n}{h}\right)^{h/2}<\left(\frac34\right)^h\frac{1}{Chn}\binom{n}{h},
\end{equation}
if $h\leq n/2$ is sufficiently large. By using the estimate 
$$\frac{1}{h}\cdot\frac{n^n}{h^h(n-h)^{n-h}}<\binom{n}{h},$$
which holds for every $10\leq h$ and $2h\leq n$ it suffices to show that
$$Ch^2ne^{3h/4}\left(\frac43\right)^h\left(\frac{n}{h}\right)^{h/2}<\left(\frac{n}{h}\right)^h\left(\frac{n}{n-h}\right)^{n-h}.$$
Setting $\gamma:=n/h\geq 2$ we shall prove that
$$Ch^2ne^{3h/4}\left(\frac43\right)^h<\gamma^{h/2}\left(\frac{\gamma}{\gamma-1}\right)^{(\gamma-1)h}.$$
Since the function $\gamma\mapsto \gamma^{1/2} \left(\frac{\gamma}{\gamma-1}\right)^{\gamma-1}$ is monotone increasing on $[2,\infty)$, the right-hand side is at least $(2\sqrt{2})^h>2.828^h$.

As $e^{3/4}\cdot \frac43<2.823$, the estimate holds if $h$ is greater than some threshold $H_1$ and $n<h^4$.

Finally, when $h^4\leq n$, the following estimate may instead be used to prove \eqref{eq-4}. As $\left(\frac{n}{h}\right)^h<\binom{n}{h}$ it suffices to show that
$$Chne^{3h/4}\left(\frac43\right)^h<\left(\frac{n}{h}\right)^{h/2},$$
which is equivalent to
$$C^{1/h}h^{1/h}e^{3/4}h^{1/2}\cdot\frac43<n^{1/2-1/h}.$$
For $h^4\leq n$ the right hand side is at least $h^{2-4/h}$, so the bound holds if $h$ is larger than some threshold $H_2$.

Hence, according to \eqref{eq-final} for every $a\in G$
$$\left|T(a)-\frac1n\binom{n}{h}\right|\leq \left(\frac34\right)^h\frac1n\cdot\frac{n!}{(n-h)!},$$
if $h>\max(10,H_1,H_2)$.

\end{proof}

\section{Concluding remarks}
We proved that $h$-element subset sums in an abelian group of order $n=2h$ are asymptotically equidistributed. Quantitatively, for every $a\in G$ there are $(1+O(0.75^h))\frac{1}{n}\binom{n}{h}$ $h$-element subsets $\{x_1,\dots,x_h\}$ where the sum of the elements is $a$. We did not attempt to optimize the constant $0.75$. The case $G=\mathbb{F}_2^r$, investigated by Katona and Makar-Limanov \cite{KML08}, shows that constant can not be decreased below 0.5.

Finally, we shall mention a paper of Wagon~and~Wilf~\cite{WW94} who looked at the question asking for characterising the cases when the $h$-element subset sums are distributed uniformly (that is, for each $a\in G$ the number of $h$-element subsets where the total sum is $a$ is {\it exactly} $\frac{1}{n}\binom{n}{h}$).

\section{Acknowledgements}
The author was supported by the Lend\"ulet program of the Hungarian Academy of Sciences (MTA) and by the National Research, Development and Innovation Office NKFIH (Grant Nr. K146387).

\medskip

\end{document}